\documentclass[letterpaper, 10 pt, conference]{ieeeconf}  

\IEEEoverridecommandlockouts                              
\usepackage{amsfonts}
\usepackage{graphicx,xcolor}
\usepackage{amsmath}
\usepackage{amssymb}
\usepackage{dsfont}

\usepackage{ifthen}
\usepackage{color}
\usepackage{cite}

\usepackage{mathtools}
\usepackage{booktabs}
\usepackage{url}
\usepackage{float}
\usepackage{multirow}

\def\scr#1{{\cal #1}}
\newcommand{\R}{{\rm I\!R}}
\def\eq#1{\begin{equation}#1\end{equation}}

\newcommand{\bbb}{\mathbb}
\newtheorem{theorem}{Theorem}
\newtheorem{lemma}{Lemma}
\newtheorem{remark}{Remark}
\newtheorem{assumption}{Assumption}
\newtheorem{definition}{Definition}
\newtheorem{proposition}{Proposition}
\newtheorem{corollary}{Corollary}

\newcommand{\dfb}{\stackrel{\Delta}{=}}
\def\qed{ \rule{.08in}{.08in}}
\DeclarePairedDelimiter\ceil{\lceil}{\rceil}
\DeclarePairedDelimiter\floor{\lfloor}{\rfloor}

\newcommand{\1}{\mathbf{1}}
\newcommand{\0}{\mathbf{0}}

\title{\LARGE \bf {\color{black} Subgradient-Push} Is of the Optimal Convergence Rate
}

\author{Yixuan Lin \hspace{.3in} Ji Liu
\thanks{
Y. Lin is with the Department of Applied Mathematics and Statistics at Stony Brook University (\texttt{yixuan.lin.1@stonybrook.edu}).
J. Liu is with the Department of Electrical and Computer Engineering at Stony Brook University
(\texttt{ji.liu@stonybrook.edu}).
}
}

\begin{document}

\maketitle
\thispagestyle{empty}
\pagestyle{empty}


\begin{abstract}
The push-sum based subgradient is an important method for distributed convex optimization over unbalanced directed graphs, which is known to converge at a rate of $O(\ln t/\sqrt{t})$. This paper shows that the {\color{black} subgradient-push} algorithm actually converges at a rate of $O(1/\sqrt{t})$, which is the same as that of the single-agent subgradient and thus optimal. The proposed tool for analyzing push-sum based algorithms is of independent interest. 
\end{abstract}


\vspace{.1in}

\section{Introduction}

There are three major information fusion schemes in the vast distributed algorithms literature: consensus via stochastic matrices \cite{reachingp1}, distributed averaging via doubly stochastic matrices \cite{fast}, and push-sum via column stochastic matrices \cite{push}.\footnote{A square nonnegative matrix is called a row stochastic matrix, or simply stochastic matrix, if its row sums all equal one. Similarly, a square nonnegative matrix is called a column stochastic matrix if its column sums all equal one.
A square nonnegative matrix is called a doubly stochastic matrix if its row sums and column sums all equal one.} Among the three, the push-sum scheme is the only one that is able to not only achieve agreement on the average, but also works for directed graphs, allowing uni-directional communication. Because of this, the push-sum scheme has been widely utilized in various distributed algorithms including distributed optimization \cite{nedic} and distributed reinforcement learning \cite{yixuan}.

The push-sum algorithm was first proposed in \cite{push} and 
sometimes also called weighted gossip \cite{weighted}, ratio consensus \cite{hadjicostis2013average}, and double linear iteration \cite{acc12}.  
Although the analysis of the push-sum algorithm is elegant, the analyses of push-sum based algorithms are often quite complicated, e.g., {\color{black} subgradient-push} \cite{nedic}, DEXTRA \cite{xi2017dextra} (a push-sum based variant of the well-known EXTRA algorithm \cite{shi2015extra}) and Push-DIGing \cite{nedic2017achieving}. Actually, all these push-sum based algorithms rely on the pioneering analysis and results in \cite{nedic}. 

Distributed optimization originated from the work of \cite{nedic2009distributed} and has
achieved great success in both theory and practice; see survey papers \cite{yang2019survey,nedic2018distributed,molzahn2017survey}. 
Most existing distributed optimization algorithms require the underlying communication network be described by an undirected graph or a balanced directed graph (a directed graph is balanced if the sum of all in-weights equals the sum of all out-weights at each of its vertices \cite{gharesifard2013distributed}), which allows a distributed manner to construct a doubly stochastic matrix. Such a distributed algorithm usually achieves the same order of convergence rate as its single-agent counterpart, with a difference at a constant coefficient depending on graph connectivity \cite{nedic2018network}.

The push-sum based subgradient algorithm proposed in \cite{nedic} is the first distributed convex optimization algorithm which works for unbalanced directed graphs. There are two ``gaps'' in the analysis in \cite{nedic}. First, the convergence rate analysis is based on a special convex combination of the history of the states of all agents (see Theorem 2 in \cite{nedic}), which is ``unusual'' compared with non-push-sum based distributed optimization algorithms (see e.g. \cite{nedic2009distributed}). Second, more importantly, the convergence rate derived in \cite{nedic} is of order $O(\ln t/\sqrt{t})$, which is slower than that of the single-agent subgradient method, $O(1/\sqrt{t})$ (see Theorem 7 in \cite{nedic2018network}). With these in mind, this paper aims to close the theoretical gap between the convergence rates of conventional single-agent subgradient and push-sum based subgradient, by analyzing the ``standard'' convex combination of the history of the states of all agents. We achieve this goal by establishing the explicit ``absolute probability sequence'' for the push-sum algorithm, which yields a novel analysis tool for push-sum based
distributed algorithms over possibly time-varying, unbalanced, directed~graphs.

\vspace{.1in}

\section{{\color{black} Subgradient-Push}}

Consider a network consisting of $n$ agents, labeled $1$ through $n$ for the purpose of presentation. The agents are not aware of such a global labeling, but can differentiate between their neighbors. The neighbor relations among the $n$ agents are characterized by a time-dependent directed graph $\bbb{G}(t) = (\mathcal{V},\mathcal{E}(t))$ whose
vertices correspond to agents and whose directed edges (or arcs) depict neighbor relations, where $\mathcal{V}=\{1,\ldots,n\}$ is the vertex set and $\mathcal{E}(t)\subset\mathcal{V} \times \mathcal{V}$ is the directed edge set at time $t$.
Specifically, agent $j$ is an in-neighbor of agent $i$ at time $t$ if $(j,i)\in\scr{E}(t)$, and similarly, agent $k$ is an out-neighbor of agent $i$ at time $t$ if $(i,k)\in\scr{E}(t)$.
Each agent can send information to its out-neighbors and receive information from its in-neighbors. Thus, the directions of edges represent the
directions of information flow. For convenience, we assume that each agent is always an in- and out-neighbor of itself, which implies that $\bbb{G}(t)$ has self-arcs at all vertices for all time $t$. We use $\mathcal{N}_i(t)$ and $\mathcal{N}_i^{-}(t)$ to denote the in- and out-neighbor set of agent $i$ at time $t$, respectively, i.e.,
\begin{align*}
    \mathcal{N}_i(t) &= \{ j \in \mathcal{V}  : ( j, i ) \in \mathcal{E}(t) \}, \\ \mathcal{N}_i^{-}(t) &= \{ k \in \mathcal{V}  :  ( i, k ) \in \mathcal{E}(t) \}.
\end{align*}
It is clear that $\mathcal{N}_i(t)$ and $\mathcal{N}_i^{-}(t)$ are nonempty as they both contain index $i$.
The goal of the $n$ agents is to cooperatively to minimize the cost function
$$f(z)=\frac{1}{n}\sum_{i=1}^n f_i(z),$$ where each $f_i$ is a ``private'' convex (not necessarily differentiable) cost function only known to agent $i$. It is assumed that the set of optimal solutions to $f$, denoted by $\scr Z$, is nonempty.

Since each $f_i$ is not necessarily differentiable, the gradient descent method may not be applicable. Instead, the subgradient method \cite{subgradient} can be applied. For a convex function $h : \R^d\rightarrow \R$, a vector $g\in\R^d$ is called a subgradient of $h$ at point $x$ if
\eq{
h(y)\ge h(x) + g^\top (y-x) \;\; {\rm for \; all} \;\; y\in\R^d.
\label{eq:subgradient}}
Such a vector $g$ always exists and may not be unique. In the case when $h$ is differentiable at point $x$, the subgradient $g$ is unique and equals $\nabla h(x)$, the gradient of $h$ at $x$. Thus, the subgradient can be viewed as a generalization of
the notion of the gradient. From \eqref{eq:subgradient} and the Cauchy-Schwarz inequality, 
\eq{
h(y) - h(x) \ge - G \| y-x\|_2,
\label{eq:G}}
where $G$ is an upper bound for the 2-norm of the subgradients of $h$ at both $x$ and $y$.

The subgradient method was first proposed in \cite{subgradient} and the first distributed subgraident method was proposed in \cite{nedic2009distributed}, which is based on average consensus.
The {\color{black} subgradient-push} algorithm, proposed in \cite{nedic},  is as follows\footnote{The algorithm is called subgradient-push in \cite{nedic} and written in a different but mathematically equivalent form there.}:
\begin{align}
    x_i(t+1) &= \sum_{j\in\scr{N}_i(t)} w_{ij}(t)\Big[x_j(t) - \alpha(t)g_j(t)\Big], \label{eq:pushsub_x}\\  y_i(t+1) &= \sum_{j\in\scr{N}_i(t)} w_{ij}(t)y_j(t),\;\;\;\;\; y_i(0)=1\label{eq:pushsub_y},
\end{align}
where $\alpha(t)$ is the stepsize, $g_j(t)$ is a subgradient of $f_j(z)$ at $x_j(t)/y_j(t)$, and $w_{ij}(t)$, $j\in\scr{N}(t)$, are positive weights satisfying the following assumption.  

\vspace{.05in}
 
\begin{assumption}\label{assum:weighted matrix}
There exists a constant $\beta>0$ such that for all $i,j\in\scr V$ and $t$, $w_{ij}(t) \ge \beta$ whenever $j\in\scr{N}_i(t)$. For all $i\in\scr V$ and $t$, $\sum_{j\in\scr{N}_i^{-}(t)} w_{ji}(t) = 1$.
\end{assumption}

\vspace{.05in}

A typical choice of $w_{ij}(t)$ is $1/|\mathcal{N}_j^{-}(t)|$ for all $j\in\scr{N}_i(t)$ which can be computed in a distributed manner and satisfies Assumption~\ref{assum:weighted matrix} with $\beta=1/n$.
Let $W(t)$ be the $n\times n$ matrix whose $ij$th entry equals $w_{ij}(t)$ if $j\in\scr{N}_i(t)$ and zero otherwise; in other words, we set $w_{ij}(t)=0$ for all $j\notin\scr{N}_i(t)$.
From Assumption~\ref{assum:weighted matrix}, each $W(t)$ is a column stochastic matrix that is compliant with the neighbor graph $\bbb{G}(t)$. 
Since each agent $i$ is always assumed to be an in-neighbor of itself, all diagonal entries of $W(t)$ are positive. 

To state the convergence result of the {\color{black} subgradient-push} algorithm, we need the following assumption and concept. 

\vspace{.05in}

\begin{assumption} \label{assum:step-size}
    The step-size sequence $\{\alpha(t)\}$ is positive, non-increasing, and satisfies $\sum_{t=0}^\infty \alpha(t) = \infty$ and $\sum_{t=0}^\infty \alpha^2(t) < \infty$. 
\end{assumption}

\vspace{.05in}

\begin{definition}\label{def:uniformly}
     A directed graph sequence $\{ \bbb{G}(t) \}$ is uniformly strongly connected if there exists a positive integer $L$ such that for any $t\ge 0$, the union graph $\cup_{k=t}^{t+L-1} \bbb{G}(k)$ is strongly connected.\footnote{A directed graph is strongly connected if it has a directed path from any vertex to any other vertex. The union of two directed graphs, $\bbb G_p$ and $\bbb G_q$, with the same vertex set, written $\bbb G_p \cup \bbb G_q$, is meant the directed graph with the same vertex set and edge set being the union of the edge set of $\bbb G_p$ and $\bbb G_q$. Since this union is a commutative and associative binary operation, the definition extends unambiguously to any finite sequence of directed graphs with the same vertex set.}
     If such an integer exists, we sometimes say that $\{ \bbb{G}(t) \}$ is uniformly strongly connected by sub-sequences of length $L$.
\end{definition}

\vspace{.05in}

It is not hard to prove that the above definition is equivalent to the two popular joint connectivity definitions in consensus literature, namely ``$B$-connected'' \cite{nedic2009distributed_quan} and ``repeatedly jointly strongly connected'' \cite{reachingp1}.

Define $z_i(t)=x_i(t)/y_i(t)$ and $\bar z(t) = \frac{1}{n}\sum_{i=1}^n z_i(t)$.

\vspace{.05in}

\begin{theorem} \label{thm:bound_everage_n_convex}
    Suppose that $\{ \bbb{G}_t \}$ is uniformly strongly connected and that $\|g_i(t)\|_2$ is uniformly bounded for all $i$ and~$t$.
\begin{itemize}
    \item[1)] If the stepsize $\alpha(t)$ is time-varying and satisfies Assumption~\ref{assum:step-size}, then 
\begin{align*}
    \lim_{t\rightarrow\infty}f\bigg(\frac{\sum_{\tau =0}^t \alpha(\tau) \bar z(\tau) }{\sum_{\tau =0}^t \alpha(\tau)}\bigg) = f(z^*).
\end{align*}
    
    \item[2)] If the stepsize is fixed and  $\alpha(t) = 1/\sqrt{T}$ for $T>0$ steps, i.e., $t\in\{0,1,\ldots,T-1\}$, then
\begin{align*}
    f\bigg(\frac{\sum_{\tau =0}^{T-1} \bar z(\tau)  }{ T } \bigg) - f(z^*) 
    & \le O\Big(\frac{1}{ \sqrt{T}}\Big).
\end{align*}
\end{itemize}

\end{theorem}

\vspace{.15in}

The above theorem establishes the convergence rate of $f((\sum_{\tau =0}^{T-1} \bar z(\tau))/T)$, as conventionally did in average-consensus-based subgradient \cite{nedic2009distributed}, and the rate is of $O(1/\sqrt{t})$, which is the same as that of the conventional single-agent subgradient method \cite[Theorem 7]{nedic2018network}. Thus, the derived convergence rate is optimal.

Theorem \ref{thm:bound_everage_n_convex} is actually a consequence of the following refined result, which further provides finite-time error bounds for the {\color{black} subgradient-push} algorithm. 

\vspace{.05in}

\begin{theorem} \label{thm:bound_everage_n_convex_bound}
    Suppose that $\{ \bbb{G}_t \}$ is uniformly strongly connected by sub-sequences of length $L$ and that $\|g_i(t)\|_2$ is uniformly bounded above by a positive number $G$ for all $i$ and $t$.
\begin{itemize}
    \item[1)] 
If the stepsize $\alpha(t)$ is time-varying and satisfies Assumption~\ref{assum:step-size}, then for all $t \ge 0$,
\begin{align}
    &\;\;\;\; f\bigg(\frac{\sum_{\tau =0}^t \alpha(\tau) \bar z(\tau) }{\sum_{\tau =0}^t \alpha(\tau)}\bigg) - f(z^*) \nonumber\\
    & \le \frac{ \| \bar z(0) - z^* \|_2^2 + G^2\sum_{\tau =0}^t \alpha^2(\tau)}{\sum_{\tau =0}^t 2\alpha(\tau) }  \nonumber\\
    &\;\;\;  + \frac{2G\alpha(0) \sum_{i=1}^n \|\bar z(0) - z_i(0) \|_2 }{n\sum_{\tau =0}^t \alpha(\tau) } \nonumber\\
    &\;\;\; {\color{black}+ \frac{32G}{\eta}  \sum_{i=1}^n \Big\| x_i(0) - \alpha(0) g_i(0) \Big\|_2 \frac{\sum_{\tau =0}^{t-1} \alpha(\tau) \mu^\tau  }{\sum_{\tau =0}^t \alpha(\tau) }} \nonumber\\
    &\;\;\; + \frac{32nG^2}{\eta(1-\mu)} \frac{\sum_{\tau =0}^{t-1} \alpha(\tau)   \left(  \alpha(0) \mu^{\tau/2}  +  \alpha(\ceil{\frac{\tau}{2}}) \right)}{\sum_{\tau =0}^t \alpha(\tau) }\label{eq:bound_timevarying}.
\end{align} 
\item[2)] If the stepsize is fixed and  $\alpha(t) = 1/\sqrt{T}$ for $T>0$ steps, i.e., $t\in\{0,1,\ldots,T-1\}$, then
\begin{align}
    &\;\;\;\; f\bigg(\frac{\sum_{\tau =0}^{T-1} \bar z(\tau)  }{ T } \bigg) - f(z^*) \nonumber\\
    & \le  \frac{2G \sum_{i=1}^n \|\bar z(0) - z_i(0) \|_2}{nT } +\frac{ \| \bar z(0) - z^* \|_2^2 + G^2 }{ 2\sqrt{T} } \nonumber\\
    &\;\;\; {\color{black}+ \frac{32 G }{\eta(1 - \mu)T}  \sum_{i=1}^n \Big\| x_i(0) - \frac{1}{\sqrt{T}} g_i(0) \Big\|_2} \nonumber\\
    &\;\;\; + \frac{32n G^2 }{ \eta(1- \mu)\sqrt{T}} \label{eq:bound_fixed}.
\end{align}
\end{itemize}
Here $\eta$ and $\mu$ are positive constants which satisfy $\eta \ge \frac{1}{n^{nL}}$ and $\mu \le (1-\frac{1}{n^{nL}})^{1/L}$, respectively, and $\ceil{\cdot}$ denotes the ceiling function.
\end{theorem}


\vspace{.05in}

The above theorem characterizes convergence rates for a network-wide averaged state. The following theorem provides convergence rates for each individual agent.

\vspace{.05in}

\begin{theorem} \label{thm:bound_everage_zi}
    Suppose that $\{ \bbb{G}_t \}$ is uniformly strongly connected by sub-sequences of length $L$ and that $\|g_i(t)\|_2$ is uniformly bounded above by a positive number $G$ for all $i$ and $t$.
\begin{itemize}
    \item[1)] 
If the stepsize $\alpha(t)$ is time-varying and satisfies Assumption~\ref{assum:step-size}, then for all $t \ge 0$ and $k\in\scr{V}$,
\begin{align}
    &\;\;\;\; f\left(\frac{\sum_{\tau =0}^t \alpha(\tau) z_k(\tau) }{\sum_{\tau =0}^t \alpha(\tau)}\right) - f(z^*) \nonumber\\
    & \le \frac{ \| \bar z(0) - z^* \|_2^2 + G^2\sum_{\tau =0}^t \alpha^2(\tau)}{\sum_{\tau =0}^t 2\alpha(\tau) }  \nonumber\\
    &\;\;  + \frac{G\alpha(0) \sum_{i=1}^n (\|\bar z(0) - z_i(0) \|_2 + \|z_k(0)  - z_i(0) \|_2)}{n\sum_{\tau =0}^t \alpha(\tau) } \nonumber\\
     &\;\; {\color{black}+ \frac{32G}{\eta}  \sum_{i=1}^n \Big\| x_i(0) - \alpha(0) g_i(0) \Big\|_2 \frac{\sum_{\tau =0}^{t-1} \alpha(\tau) \mu^\tau  }{\sum_{\tau =0}^t \alpha(\tau) }} \nonumber\\
     &\;\; + \frac{32nG^2}{\eta(1-\mu)} \frac{\sum_{\tau =0}^{t-1} \alpha(\tau)   \left(  \alpha(0) \mu^{\tau/2}  +  \alpha(\ceil{\frac{\tau}{2}}) \right)}{\sum_{\tau =0}^t \alpha(\tau) }\label{eq:bound_timevarying_zi}.
\end{align} 
\item[2)] If the stepsize is fixed and  $\alpha(t) = 1/\sqrt{T}$ for $T>0$ steps, i.e., $t\in\{0,1,\ldots,T-1\}$, then for any $k\in\scr{V}$,
\begin{align}
    &\;\;\;\; f\left(\frac{\sum_{\tau =0}^{T-1} z_k(\tau)  }{ T } \right) - f(z^*) \nonumber\\
    & \le  \frac{ \| \bar z(0) - z^* \|_2^2 + G^2 }{ 2\sqrt{T} }+ \frac{32n G^2 }{ \eta(1- \mu)\sqrt{T}} \nonumber\\
    &\;\;\; + \frac{G \sum_{i=1}^n \|\bar z(0)  - z_i(0) \|_2+\|z_k(0) - z_i(0) \|_2}{nT } \nonumber\\
    &\;\;\; {\color{black}+ \frac{32 G }{\eta(1 - \mu)T}  \sum_{i=1}^n \Big\| x_i(0) - \frac{1}{\sqrt{T}} g_i(0) \Big\|_2 } \label{eq:bound_fixed_zi}.
\end{align}
\end{itemize}
Here the positive constants $\eta$ and $\mu<1$ are the same as in Theorem \ref{thm:bound_everage_n_convex_bound}.
\end{theorem}

\vspace{.05in}

Using the same argument as in the proof
of Theorem \ref{thm:bound_everage_n_convex}, we have for each agent $k\in\scr V$, with a time-varying stepsize $\alpha(t)$ satisfying Assumption~\ref{assum:step-size},
\begin{align*}
    \lim_{t\rightarrow\infty}f\bigg(\frac{\sum_{\tau =0}^t \alpha(\tau) z_k(\tau) }{\sum_{\tau =0}^t \alpha(\tau)}\bigg) = f(z^*),
\end{align*}
and with a fixed stepsize $\alpha(t) = 1/\sqrt{T}$ for $T>0$ steps, 
\begin{align*}
    f\bigg(\frac{\sum_{\tau =0}^{T-1} z_k(\tau)  }{ T } \bigg) - f(z^*) 
    & \le O\Big(\frac{1}{ \sqrt{T}}\Big).
    \end{align*}

\vspace{.1in}

\section{Analysis}

In this section, we provide a novel analysis of the {\color{black} subgradient-push} algorithm \eqref{eq:pushsub_x}--\eqref{eq:pushsub_y} and proofs of Theorems \ref{thm:bound_everage_n_convex} and~\ref{thm:bound_everage_n_convex_bound}. The analysis appeals to the concept of ``absolute probability sequence'' for push-sum. Thus, we begin with revisiting the well-known push-sum algorithm.

\subsection{Push-Sum}

In the push-sum algorithm, each agent $i$ has control over two variables, $x_i(t)$ and $y_i(t)$, which are updated as follows:  
\begin{align}
    x_i(t+1) &= \sum_{j\in\scr{N}_i(t)} w_{ij}(t)x_j(t),\label{pushsumx}\\ 
    y_i(t+1) &= \sum_{j\in\scr{N}_i(t)} w_{ij}(t)y_j(t),\;\;\;\;\; y_i(0)=1, \label{pushsumy}
\end{align}
where $w_{ij}(t)$, $j\in\scr{N}(t)$, are positive weights satisfying Assumption \ref{assum:weighted matrix}.

{\color{black}  Let $x(t) \dfb [x_1(t)\;\cdots\; x_n(t)]^\top\in\R^{n\times d}$ and $y(t)$ be the vector in $\R^n$ whose $i$th entry is  $y_i(t)$. From \eqref{pushsumx} and \eqref{pushsumy}, $x(t+1)=W(t)x(t)$ and $y(t+1)=W(t)y(t)$. Since $W(t)$ is always column stochastic for all $t\ge 0$, it is easy to show that $\sum_{i=1}^n x_i(t)=\sum_{i=1}^n x_i(0)$ and $\sum_{i=1}^n y_i(t)=\sum_{i=1}^n y_i(0)=n$ for all $t\ge 0$.} 

\vspace{.05in}

\begin{lemma} \label{lemma:pushsum_product}
    Suppose that $\{ \bbb{G}(t) \}$ is uniformly strongly connected. Then, for any fixed $\tau\ge 0$, 
    $W(t)\cdots W(\tau+1)W(\tau)$
    will converge to the set 
    $\{v\1^\top \; :\; v\in\R^n, \1^\top v=1, v>\0\}$ exponentially fast as $t\rightarrow\infty$.\footnote{We use $\0$ and $\1$ to denote the vectors whose entries all equal to $0$ or $1$, respectively, where the dimensions of the vectors are to be understood from the context. We use $v>\0$ to denote a positive vector, i.e., each entry of $v$ is positive.}
\end{lemma}

\vspace{.05in}

The lemma is essentially the same as Corollary~2~(a) in~\cite{nedic}. Suppose $\{ \bbb{G}_t \}$ is uniformly strongly connected by sub-sequences of length $L$, Lemma \ref{lemma:pushsum_product}  implies that there exist constants $c>0$ and $\mu \in [0,1)$ and a sequence of stochastic vectors\footnote{A vector is called a stochastic vector if
its entries are all nonnegative and sum to one.} $ \{ v(t)\}$ such that for all $i,j \in \mathcal{V}$ and $t \ge \tau \ge 0$,
\begin{align}\label{mu}
    \big| \big[W(t)\cdots W(\tau+1)W(\tau)\big]_{ij} - {\color{black}v_i(t)} \big|\le c \mu^{t-\tau},
\end{align}
where $[\cdot]_{ij}$ denotes the $ij$th entry of a matrix. 
In \cite{nedic}, it has been shown that $c=4$ and $\mu=(1-\frac{1}{n^{nL}})^{1/L}$.




To proceed, we define a time-dependent $n\times n$ matrix $S(t)$ whose $ij$th entry is 
\begin{equation}\label{eq:s}
    s_{ij}(t)= \frac{w_{ij}(t)y_j(t)}{y_i(t+1)}=\frac{w_{ij}(t)y_j(t)}{\sum_{k=1}^n w_{ik}(t)y_k(t)}.
\end{equation}
It is worth emphasizing that $S(t)$ is independent of $x(t)$.
The following lemma guarantees that $S(t)$ is well defined.

\vspace{.05in}

\begin{lemma}\label{lemma:y_bound}
    Suppose that $\{ \bbb{G}(t) \}$ is uniformly strongly connected, then there exists a constant $\eta>0$ such that $n \ge y_i(t) \ge \eta$ for all $i$ and $t$.
\end{lemma}

\vspace{.05in}

The lemma is essentially the same as Corollary~2~(b) in~\cite{nedic}, which further proves that if $\{ \bbb{G}_t \}$ is uniformly strongly connected by sub-sequences of length $L$, then 
$\eta \ge \frac{1}{n^{nL}}$.



Define $z_i(t)=x_i(t)/y_i(t)$ for each $i\in\mathcal{V}$. Then, 
\begin{align}
&\;\;\;\;\; z_i(t+1) 
= \frac{x_i(t+1)}{y_i(t+1)} 
= \frac{\sum_{j=1}^n w_{ij}(t)x_j(t)}{\sum_{j=1}^n w_{ij}(t)y_j(t)}\nonumber\\
&= \sum_{j=1}^n \frac{w_{ij}(t)x_j(t)}{\sum_{k=1}^n w_{ik}(t)y_k(t)}
= \sum_{j=1}^n \bigg[\frac{w_{ij}(t)y_j(t)}{\sum_{k=1}^n w_{ik}(t)y_k(t)}\bigg]z_j(t) \nonumber\\
&= \sum_{j=1}^n s_{ij}(t)z_j(t),\label{eq:update_ratio}
\end{align}
which implies that $z(t+1)=S(t)z(t)$ where $z(t)$ is the vector in $\R^n$ whose $i$th entry is $z_i(t)$.
Actually $S(t)$ is always a stochastic matrix, as we will show shortly. 

Similar to the discrete-time state transition matrix, 
let $\Phi_W(t,\tau)=W(t-1)\cdots W(\tau)$ with $t>\tau$,
and similarly, let $\Phi_S(t,\tau)=S(t-1)\cdots S(\tau)$ with $t>\tau$.

\vspace{.05in}

\begin{lemma}\label{lemma:yixuan}
For $i,j\in\mathcal{V}$ and $t > \tau \ge 0$, there holds
$[\Phi_S(t,\tau)]_{ij}y_i(t)=[\Phi_W(t,\tau)]_{ij}y_j(\tau).$
\end{lemma}

\vspace{.05in}

{\em Proof of Lemma~\ref{lemma:yixuan}:}
The claim will be proved by induction on $t$. 
For the basis step, suppose that $t=\tau+1$. Then, from \eqref{eq:s},
$[\Phi_S(\tau+1,\tau) ]_{ij} = s_{ij}(\tau) = \frac{ y_j(\tau)  w_{ij}(\tau)}{y_i(\tau+1)} = \frac{ y_j(\tau)}{y_i(\tau+1)} [\Phi_W(\tau+1,\tau) ]_{ij}.
$
Thus, in this case the claim is true.
For the inductive step, 
suppose that the claim holds for $t=h > \tau$, where $h$ is a positive integer, and that $t=h+1$. Then,
\begin{align*}
&\;\;\;\;\;[\Phi_S(h+1,\tau)]_{ij} 
= \sum_{k=1}^n  s_{ik}(h) \cdot [\Phi_S(h,\tau)]_{kj} \\
&= \sum_{k=1}^n  \frac{w_{ik}(h)y_k(h)}{y_i(h+1)}  \cdot \frac{ y_j(\tau)  }{y_k(h)} [\Phi_W(h,\tau) ]_{kj}  \\
&= \frac{ y_j(\tau)  }{y_i(h+1)} \sum_{k=1}^n w_{ik}(h) \cdot  [\Phi_W(h,\tau) ]_{kj} \\
&= \frac{ y_j(\tau)  }{y_i(h+1)} [\Phi_W(h+1,\tau) ]_{ij} , 
\end{align*}
which establishes the claim by induction.
\hfill$\qed$

\vspace{.05in}

More can be said.

\vspace{.05in}

\begin{lemma}\label{yixuan}
    Suppose that $\{ \bbb{G}(t) \}$ is uniformly strongly connected. Then, 
    for any fixed $\tau\ge0$, $S(t)\cdots S(\tau+1)S(\tau)$ will converge to $\frac{1}{n} \1 y^\top(\tau)$. 
\end{lemma}

\vspace{.05in}

{\em Proof of Lemma~\ref{yixuan}:}
From Lemma~\ref{lemma:pushsum_product}, for any given $\tau\ge 0$, there holds $\lim_{t\to\infty} [\Phi_W(t,\tau)] = v(\tau,\infty) \1^\top $, with the understanding that $v(\tau,\infty)$ is not necessarily a constant vector. From Lemma~\ref{lemma:yixuan} and the fact that $y(t) =  \Phi_W(t,\tau) y(\tau)$ for all $t > \tau$, for any $i, j\in \mathcal{V}$ we have
\begin{align*}
&\;\;\;\;\; \lim_{t\to\infty} [\Phi_S(t,\tau)]_{ij}  \\
&= \lim_{t\to\infty} \frac{ y_j(\tau)  }{y_i(t)} [\Phi_W(t,\tau) ]_{ij} 
= \lim_{t\to\infty} \frac{ y_j(\tau) [\Phi_W(t,\tau) ]_{ij} }{\sum_{k=1}^n [\Phi_W(t,\tau) ]_{ik} y_k(\tau)} \\
&=  \frac{ y_j(\tau) \lim_{t\to\infty} [\Phi_W(t,\tau) ]_{ij} }{\lim_{t\to\infty} \sum_{k=1}^n [\Phi_W(t,\tau) ]_{ik} y_k(\tau)} 
=  \frac{ y_j(\tau) v_i(\tau,\infty) }{  \sum_{k=1}^n v_i(\tau,\infty) y_k(\tau)}\\ 
&\overset{(a)}{=}  \frac{ y_j(\tau) }{  \sum_{k=1}^n y_k(\tau)} 
\overset{(b)}{=} \frac{y_j(\tau)}{n},
\end{align*}
where in (a) we used the fact that $v(\tau,\infty)>\0$ by Lemma~\ref{lemma:pushsum_product} and in (b) we used the fact that $\sum_{i=1}^n y_i(t) = n $ for all $t \ge 0$.~\hfill$\qed$

\vspace{.05in}

\begin{proposition} \label{lemma:expconvergen} 
    Suppose that $\{ \bbb{G}(t) \}$ is uniformly strongly connected. Then, for any fixed $\tau\ge0$, $S(t)\cdots S(\tau+1)S(\tau)$ will converge to $\frac{1}{n} \1 y^\top(\tau)$ exponentially fast as $t\rightarrow\infty$. 
\end{proposition}

\vspace{.05in}

{\em Proof of Proposition~\ref{lemma:expconvergen}:}
From \eqref{mu}, there exist constants $c>0$ and $\mu \in [0,1)$ and a sequence of stochastic vectors $ \{ v(t)\}$ such that 
$
    | [\Phi_W(t+1,\tau)]_{ij} - v(t) |\le c \mu^{t-\tau}
$ for all $i,j \in \mathcal{V}$ and $t \ge \tau \ge 0$.
Recall that $\sum_{i=1}^n y_i(t)$ always equals $n$ and, by Lemma~\ref{lemma:y_bound},  all $y_i(t)$ are always positive. From Lemma~\ref{lemma:yixuan}, for all $t\ge\tau\ge0$,
\begin{align*}
&\;\;\;\;\;\Big| [\Phi_S(t+1,\tau)]_{ij} - \frac{y_j(\tau)}{n} \Big| \\
&= \Big| \frac{ y_j(\tau) [\Phi_W(t+1,\tau) ]_{ij} }{y_i(t+1)}  - \frac{y_j(\tau)}{n} \Big|\\
&=\Big| \frac{ ny_j(\tau) [\Phi_W(t+1,\tau) ]_{ij} - y_j(\tau) [\Phi_W(t+1,\tau)y(\tau)]_i}{ny_i(t+1)}   \Big|
\\
&= \Big| \frac{n y_j(\tau) \left([\Phi_W(t+1,\tau)]_{ij} - v_i(t) + v_i(t) \right) }{ny_i(t+1)} \\
&\;\;\;\; -\frac{ y_j(\tau)\sum_{k=1}^n ([\Phi_W(t+1,\tau)]_{ik} - v_i(t) + v_i(t) ) y_k(\tau) }{ny_i(t+1)}\Big|,
\end{align*}
which implies that
\begin{align*}
&\;\;\;\;\;\Big| [\Phi_S(t+1,\tau)]_{ij} - \frac{y_j(\tau)}{n} \Big| \\
&=  \Big| \frac{ n y_j(\tau) ([\Phi_W(t+1,\tau)]_{ij} - v_i(t) )    }{ ny_i(t+1) } \\
&\;\;\;\; -\frac{ y_j(\tau) \sum_{k=1}^n ( [\Phi_W(t+1,\tau)]_{ik}- v_i(t) ) y_k(\tau)    }{ ny_i(t+1) } \Big|\\
&\le    \frac{  n y_j(\tau) \left|[\Phi_W(t+1,\tau)]_{ij} - v_i(t) \right| }{ ny_i(t+1) } \\
&\;\;\;\; + \frac{   y_j(\tau) \sum_{k=1}^n \left| [\Phi_W(t+1,\tau)]_{ik}- v_i(t) \right| y_k(\tau) }{ ny_i(t+1) }\\
& \le  \frac{ n y_j(\tau) c \mu^{t-\tau} + y_j(\tau) \sum_{k=1}^n c \mu^{t-\tau} y_k(\tau)     }{ ny_i(t+1) } \\
&= \frac{ 2 y_j(\tau) c \mu^{t-\tau} }{ y_i(t+1) } \le  \frac{2cn}{ \eta } \mu^{t-\tau},
\end{align*}
where we used Lemma~\ref{lemma:y_bound} in the last inequality. The above immediately implies the proposition.
\hfill$\qed$

\vspace{.05in}

The proposition immediately implies the following results.

\vspace{.05in}

\begin{corollary}\label{coro:sproduct}
    Suppose that $\{ \bbb{G}(t) \}$ is uniformly strongly connected. Then, $S(t)\cdots S(1)S(0)$ will converge to $\frac{1}{n}\1\1^\top$ exponentially fast as $t\rightarrow\infty$.
\end{corollary}

\vspace{.05in}

{\em Proof of Corollary \ref{coro:sproduct}:} The corollary is a special case of Proposition~\ref{lemma:expconvergen} by setting $\tau=0$. 
\hfill$\qed$

\vspace{.05in}

\begin{corollary}\label{thm:pushsum}
    If $\{ \bbb{G}(t) \}$ is uniformly strongly connected, then $x_i(t)/y_i(t)$ for all $i\in\scr V$ converges to $\frac{1}{n}\sum_{i=1}^n x_i(0)$ exponentially fast. 
\end{corollary}

\vspace{.05in}

{\em Proof of Corollary \ref{thm:pushsum}:}
From \eqref{eq:update_ratio}, $z(t+1)=S(t)z(t)=S(t)\cdots S(0)z(0)=S(t)\cdots S(0)x(0)$. From Corollary~\ref{coro:sproduct}, $z(t+1)$ will converge to $\frac{1}{n}\1\1^\top x(0) = (\frac{1}{n}\sum_{i=1}^nx_i(0))\1$ exponentially fast as $t\rightarrow\infty$, which completes the proof.
\hfill$\qed$

\vspace{.05in}

Although the above proof of Corollary~\ref{thm:pushsum} looks more complicated than the conventional convergence proof of the push-sum algorithm (e.g., \cite{push,hadjicostis2013average,acc12}), it yields the following novel and key property of push-sum.

To proceed, we rewrite the push-sum algorithm in a different form which directly characterizes the dynamics of $z_i(t)=x_i(t)/y_i(t)$. 
From \eqref{eq:s} and \eqref{eq:update_ratio}, $z_i(t+1)=\sum_{j=1}^n s_{ij}(t)z_j(t)$ and $s_{ij}(t)$ satisfies the following assumption.

\begin{assumption}\label{assum:rowstochastic}
     There exists a constant $\gamma>0$ such that for all $i,j\in \mathcal{V}$ and $t$, $s_{ii}(t) \ge \gamma$ and $s_{ij}(t) \ge \gamma$ whenever $s_{ij}(t)>0$. For all $i\in \mathcal{V}$ and $t$, $\sum_{j=1}^n s_{ij}(t) = 1$.
\end{assumption}

\vspace{.05in}

\begin{lemma}\label{lemma:smatrices}
    Suppose that Assumption~\ref{assum:weighted matrix} holds. Then, 
    $s_{ij}(t)$ satisfies Assumption~\ref{assum:rowstochastic} for each $t\ge 0$.
\end{lemma}

\vspace{.05in}

{\em Proof of Lemma~\ref{lemma:smatrices}:}
From Assumption~\ref{assum:weighted matrix}, each $ W(t)$ is a column stochastic matrix whose diagonal entries are all positive and $ w_{ij}(t) \ge \beta$ whenever $ w_{ij}(t) > 0$. From \eqref{eq:s}, $s_{ij}(t)>0$ only if $w_{ij}(t)>0$. From Lemma~\ref{lemma:y_bound}, when $w_{ij}(t)>0$, 
{\color{black}
$$s_{ij}(t)=\frac{w_{ij}(t)y_j(t)}{y_i(t+1)}
\ge \frac{ \beta \eta}{n}.$$ }
The above inequality and Assumption~\ref{assum:weighted matrix} imply that $s_{ij}(t)$ satisfies the first sentence of Assumption~\ref{assum:rowstochastic} with $\gamma=\beta\eta/n$.
For the second sentence of Assumption~\ref{assum:rowstochastic}, it is easy to see that 
\begin{align*}
    \sum_{j=1}^n s_{ij}(t) & = \sum_{j=1}^n \frac{w_{ij}(t)y_j(t)}{\sum_{k=1}^n w_{ik}(t)y_k(t)}  = 1
\end{align*}
for all $i\in \mathcal{V}$ and $t$, which completes the proof.
\hfill$\qed$

\vspace{.05in}

From Lemma~\ref{lemma:smatrices}, each $S(t)$ is a row stochastic matrix whose diagonal entries are all positive and whose nonzero entries are all uniformly bounded below by some positive number. More can be said. The following lemma shows that each $S(t)$ is compliant with the neighbor graph $\bbb{G}(t)$. 

\vspace{.05in}

\begin{lemma}\label{lemma:sgraph}
    The graph of $S(t)$ is the same as the graph of $W(t)$ for all $t$.\footnote{The graph of an $n\times n$ matrix is a direct graph with $n$ vertices and an arc from vertex $i$ to vertex $j$ whenever the $ji$th entry of the matrix is~nonzero.}
\end{lemma}

\vspace{.05in}

{\em Proof of Lemma~\ref{lemma:sgraph}:}
From \eqref{eq:s} and  Lemma~\ref{lemma:y_bound}, it is easy to see that $s_{ij}(t)>0$ if and only if $w_{ij}(t)$, which proves the lemma.
\hfill$\qed$

\vspace{.05in}

From~\eqref{eq:update_ratio}, $z(t+1)=S(t)z(t)$. 
The above lemmas imply that the dynamics of $z(t)$ is a nonlinear consensus process as $S(t)$ is dependent on $z(t)$. 
Such a transition in analysis from $x(t)$ dynamics to $z(t)$ dynamics has been used in \cite{iutzeler2013analysis}. 
To analyze such a process, we appeal to the following concept. To our knowledge, the concept has never been used to analyze the push-sum algorithm and its applications.  

\vspace{.05in}

\begin{definition}\label{def: absolute prob}
Let $\{ S(t) \}$ be a sequence of stochastic matrices. A sequence of stochastic vectors $\{ \pi(t) \}$ is an absolute probability sequence for $\{ S(t) \}$ if
$\pi^\top(t) = \pi^\top(t+1) S(t)$ for all $t\ge0$.
\end{definition}

\vspace{.05in}

This definition was first introduced by Kolmogorov \cite{kolmogorov}. 
It was shown by Blackwell \cite{blackwell} that every sequence of stochastic matrices
has an absolute probability sequence.
In general, a sequence of stochastic matrices may have more than one absolute probability sequence; when the sequence of stochastic matrices is ``ergodic'',\footnote{
A sequence of stochastic matrices $\{S(t)\}$ is called ergodic if $\lim_{t\rightarrow\infty}S(t)\cdots S(\tau+1)S(\tau)$ exists for all $\tau$.} it has a unique absolute probability sequence \cite[Lemma~1]{tacrate}. It is easy to see that when $S(t)$ is a fixed irreducible stochastic matrix $S$, $\pi(t)$ is simply the normalized left eigenvector of $S$ for eigenvalue one, and when $\{S(t)\}$ is an ergodic sequence of doubly stochastic matrices, $\pi(t)=(1/n)\1$.
More can be said. 

\vspace{.05in}

\begin{lemma} \label{lemma:bound_pi_jointly}
    {\rm (Theorem 4.8 in \cite{touri2012product})}
    Let $\{S(t)\}$ be a sequence of stochastic matrices satisfying Assumption~\ref{assum:rowstochastic}. If the graph sequence of $\{\bbb{G}(t)\}$ is uniformly strongly connected, then there exists a unique absolute probability sequence $\{ \pi(t) \}$ for the matrix sequence $\{S(t)\}$ and a constant $\pi_{\min} \in (0,1)$ such that $\pi_i(t) \ge \pi_{\min}$ for all $i$ and $t$.
\end{lemma}

\vspace{.05in}

A particular important property of the absolute probability sequence for $\{S(t)\}$ is as follows. 

\vspace{.05in}

\begin{proposition} \label{lemma:push-sum_pi_intfty}
    Suppose that $\{ \bbb{G}(t) \}$ is uniformly strongly connected. Then, the sequence of stochastic matrices $\{S(t)\}$ has a unique absolute probability sequence $\{\pi(t)\}$ with
    $\pi_i(t)=\frac{y_i(t)}{n}$ for all $i\in \mathcal{V}$ and $t \ge 0$. 
\end{proposition}

\vspace{.05in}



The proposition is a consequence of Lemma 1 in \cite{tacrate}. We provide two alternative proofs. 

\vspace{.05in}

{\em Proof of Proposition~\ref{lemma:push-sum_pi_intfty}:}
First, Lemma~\ref{yixuan} shows that $\{S(t)\}$ is ergodic, so it must have a unique absolute probability sequence $\{\pi(t)\}$. 
From Definition~\ref{def: absolute prob} and Lemma~\ref{yixuan}, for any $\tau \ge 0$,
\begin{align*}
\pi^\top(\tau) 
&= \pi^\top(\tau+1) S(\tau) 
= \pi^\top(\tau+2) S(\tau+1)S(\tau)\\
&= \lim_{t\to\infty}\pi^\top(t+1) S(t)\cdots S(\tau+1)S(\tau) \\
&= \lim_{t\to\infty} \frac{1}{n}\pi^\top(t+1)  \1 y^\top(\tau) 
= \frac{1}{n} y^\top(\tau),
\end{align*}
which proves the statement. 

Alternatively, we can also prove the proposition by showing that the sequence $\{\pi(t)\}$ with 
$\pi_i(t) = \frac{y_i(t)}{n}$ satisfies $\pi^\top(t) = \pi^\top(t+1) S(t)$. To see this, from \eqref{eq:s} and Assumption~\ref{assum:weighted matrix}, for $j \in \mathcal{V}$ 
\begin{align*}
    &\;\;\;\; [\pi^\top(t+1) S(t)]_j 
    = \sum_{i=1}^n \frac{y_i(t+1)}{n} s_{ij}(t)  \\
    &= \sum_{i=1}^n \frac{y_i(t+1)}{n} \frac{w_{ij}(t)y_j(t)}{y_i(t+1)} 
    = \sum_{i=1}^n  \frac{w_{ij}(t)y_j(t)}{n} \\
    &= \frac{y_j(t)}{n} = \pi_j^\top(t).
\end{align*}
This completes the proof.
\hfill$\qed$

\vspace{.05in}

Next we will appeal to this property to construct a novel time-varying Lyapunov function for distributed convex optimization which yields an improved convergence rate of the {\color{black} subgradient-push} algorithm.

\begin{remark}
Since the stochastic matrix sequence $S(t)$ defined by \eqref{eq:s} is purely based on the $y_i(t)$ variables and is thus
independent of the $x_i(t)$ variables of the push-sum algorithm, so its absolute probability sequence. Considering the fact that the push-sum and {\color{black} subgradient-push} algorithms share the same $y_i(t)$ dynamics which is independent of their $x_i(t)$ dynamics, all the results of $\{S(t)\}$ and its absolute probability sequence derived in this subsection also apply to the {\color{black} subgradient-push} algorithm. 
\hfill$\Box$
\end{remark}

\subsection{{\color{black} Subgradient-Push}}

We first rewrite the {\color{black} subgradient-push} algorithm as follows. 
From \eqref{eq:pushsub_x}--\eqref{eq:pushsub_y}, we have
\begin{align*}
z_i(t+1) 
&= \frac{x_i(t+1)}{y_i(t+1)} 
= \frac{\sum_{j=1}^n w_{ij}(t)[x_j(t)- \alpha(t)g_j(t)]}{\sum_{j=1}^n w_{ij}(t)y_j(t)}\\
&= \sum_{j=1}^n \frac{w_{ij}(t)[x_j(t)- \alpha(t)g_j(t)]}{\sum_{k=1}^n w_{ik}(t)y_k(t)} \\
&= \sum_{j=1}^n \bigg[\frac{w_{ij}(t)y_j(t)}{\sum_{k=1}^n w_{ik}(t)y_k(t)}\bigg]\bigg[z_j(t)- \alpha(t)
\frac{g_j(t)}{y_j(t)}\bigg]\nonumber\\
&= \sum_{j=1}^n s_{ij}(t)\bigg[z_j(t)- \alpha(t)
\frac{g_j(t)}{y_j(t)}\bigg],
\end{align*}
where $s_{ij}(t)$ is defined in \eqref{eq:s}.
In addition, 
\begin{align*}
    \bar z(t+1) 
    &= \frac{1}{n}\sum_{i=1}^n z_i(t+1) \\
    &= \frac{1}{n}\sum_{i=1}^n \sum_{j=1}^n s_{ij}(t)\bigg[z_j(t)- \alpha(t) \frac{g_j(t)}{y_j(t)}\bigg].
\end{align*}
Define a time-varying Lyapunov function 
$$\langle z(t) \rangle= \pi^\top(t)z(t).$$ 
Then, from Definition~\ref{def: absolute prob}, we have
\begin{align*}
    &\;\;\;\;\; \langle z(t+1) \rangle 
    = \sum_{i=1}^n \pi_i({t+1}) z_i(t+1) \\
    &= \sum_{i=1}^n \sum_{j=1}^n \pi_i({t+1})  s_{ij}(t)\bigg[z_j(t)- \alpha(t) \frac{g_j(t)}{y_j(t)}\bigg]\\
    &= \sum_{j=1}^n \pi_j(t) \bigg[z_j(t)- \alpha(t) \frac{g_j(t)}{y_j(t)}\bigg]=\langle z(t) \rangle - \frac{\alpha(t)}{n}\sum_{i=1}^n g_i(t),
\end{align*}
where we use the Proposition~\ref{lemma:push-sum_pi_intfty} in the last equality.

To prove Theorem~\ref{thm:bound_everage_n_convex}, we need the following lemma.

\vspace{.05in}

 \begin{lemma} \label{lemma:bound_consensus_push_SA} 
    Suppose that $\{ \bbb{G}_t \}$ is uniformly strongly connected by sub-sequences of length $L$ and that $\|g_i(t)\|_2$ is uniformly bounded above by a positive number $G$ for all $i$ and $t$.
    Then, for all $t \ge 0$ and $i \in \mathcal{V}$,
 \begin{align*}
    &\;\;\;\;\;\| z_i(t+1) - \frac{1}{n} \sum_{j=1}^n (x_j(t) -\alpha(t) g_j(t))\|_2 \\
    &{\color{black}\le \frac{8}{\eta}  \mu^t  \sum_{i=1}^n \| x_i(0) - \alpha(0) g_i(0) \|_2 +  \frac{8nG}{\eta}  \sum_{s=1}^t \mu^{t-s} \alpha(s).}
\end{align*}
Suppose, in addition, that Assumption~\ref{assum:step-size} holds. Then, for all $t \ge 0$ and $i \in \mathcal{V}$,
\begin{align*}
    &\;\;\;\;\;\| z_i(t+1) - \frac{1}{n} \sum_{j=1}^n (x_j(t) -\alpha(t) g_j(t))\|_2 \\
    & {\color{black}\le \frac{8}{\eta}  \mu^t  \sum_{i=1}^n \| x_i(0) - \alpha(0) g_i(0) \|_2 }\\
    & \;\;\; +  \frac{8nG}{\eta(1-\mu)}   \left(  \alpha(0) \mu^{t/2}  +  \alpha(\ceil{\frac{t}{2}}) \right).
\end{align*}
Here $\eta>0$ and $\mu\in(0,1)$ are defined in Lemma~\ref{lemma:y_bound} and \eqref{mu}, respectively.
\end{lemma}

\vspace{.05in}

{\em Proof of Lemma~\ref{lemma:bound_consensus_push_SA}:}
Let $h(t) = x(t) -\alpha(t) g(t) $. 
Then, 
\begin{align*}
    h(t+1) 
    &= x(t+1) -\alpha(t+1) g(t+1) \\
    &= W(t) h(t) -\alpha(t+1) g(t+1) \\
    & = \Phi_W(t+1,0) h(0) - \sum_{l=1}^t \alpha(l) \Phi_W(t+1,l) g(l) \\
    &\;\;\;\;-\alpha(t+1) g(t+1).
\end{align*}
In addition, 
\begin{align}
    &\;\;\;\;\;W(t+1) h(t+1)\nonumber \\
    & = \Phi_W(t+2,0) h(0) - \sum_{l=1}^{t+1} \alpha(l) \Phi_W(t+2,l) g(l)\label{eq:h_update_tplus1}
\end{align}
and 
\begin{align}
    \1^\top h(t+1) 
    & = \1^\top h(0) - \sum_{l=1}^{t+1} \alpha(l) \1^\top g(l) \label{eq:h_update_1}.
\end{align}
From Lemma~\ref{lemma:pushsum_product} and \eqref{mu}, there exists a sequence $ \{ \phi(t)\}$ of stochastic vectors, such that for all $i,j \in \mathcal{V}$ and $t\ge s \ge 0$
\begin{align*}
    | [\Phi_W(t+1,s)]_{ij} - \phi_i(t) |\le 4 \mu^{t-s}.
\end{align*}
Let $ D(s:t) = \Phi_W(t+1,s) - \phi(t) \1^\top$. 
From \eqref{eq:h_update_tplus1} and \eqref{eq:h_update_1}, 
\begin{align*}
    &\;\;\;\; W(t+1) h(t+1)  - \phi(t+1)\1^\top h(t+1)\\
    & = \Phi_W(t+2,0) h(0) - \sum_{l=1}^{t+1} \alpha(l) \Phi_W(t+2,l) g(l)\\
    &\;\;\;\; - \phi(t+1) (\1^\top h(0) - \sum_{l=1}^{t+1}  \alpha(l) \1^\top g(l)) \\
    & = (\Phi_W(t+2,0) - \phi(t+1) \1^\top) h(0) \\
    &\;\;\;\;- \sum_{l=1}^{t+1} \alpha(l)( \Phi_W(t+2,l) - \phi(t+1) \1^\top) g(l) \\
    & = D(0:t+1) h(0) - \sum_{l=1}^{t+1} \alpha(l) D(l:t+1) g(l),
\end{align*}
which implies that
\begin{align*}
    &\;\;\;\;\;x(t+1) = W(t) h(t) \\
    &= \phi(t)\1^\top h(t) + D(0:t) h(0) - \sum_{l=1}^{t} \alpha(l) D(l:t) g(l).
\end{align*}
Moreover,
$
    y(t+1) = \Phi_W(t+1,0) y(0) = D(0:t) \1 + n \phi(t).
$ 
Thus, for all $i\in\scr V$
\begin{align*}
    & \;\;\;\;\; z_i(t+1) - \frac{ h(t)^\top \1}{n} = \frac{x_i(t+1)}{y_i(t+1)} - \frac{\1^\top h(t)}{n}\\
    &= \frac{\phi_i(t) h(t)^\top \1 + \sum_{k=1}^n [D(0:t)]_{ik} h_k(0) }{[D(0:t) \1]_i + n \phi_i(t)}\\
    &\;\;\;\; -\frac{ \sum_{l=1}^{t} \alpha(l) \sum_{k=1}^n [D(l:t+1)]_{ik} g_k(l)}{[D(0:t) \1]_i + n \phi_i(t)} - \frac{ h(t)^\top \1}{n} \\
    &= \frac{n \sum_{k=1}^n [D(0:t)]_{ik} h_k(0)   - [D(0:t) \1]_i h(t)^\top \1}{n[D(0:t) \1]_i + n^2 \phi_i(t)}\\
    &\;\;\;\; -\frac{  n\sum_{l=1}^{t} \alpha(l) \sum_{k=1}^n [D(l:t+1)]_{ik} g_k(l)   }{n[D(0:t) \1]_i + n^2 \phi_i(t)}.
\end{align*}
From the definition of $D(0:t)$, we have $[D(0:t) \1]_i + n \phi_i(t) = [\Phi_W(t+1,0)\1]_i \ge \eta$. Therefore, 
\begin{align*}
    &\;\;\;\; \| z_i(t+1) - \frac{ h(t)^\top \1}{n} \|_2 \\
    & \le \frac{n \|\sum_{k=1}^n [D(0:t)]_{ik} h_k(0)\|_2  + \|[D(0:t) \1]_i h(t)^\top \1 \|_2 }{n[D(0:t) \1]_i + n^2 \phi_i(t)}\\
    & \;\;\;\; + \frac{n\sum_{l=1}^{t} \alpha(l) \|\sum_{k=1}^n [D(l:t)]_{ik} g_k(l)\|_2  }{n[D(0:t) \1]_i + n^2 \phi_i(t)}\\
    &{\color{black} \le \frac{n  (\max_k [D(0:t)]_{ik}) \sum_{k=1}^n  \| h_k(0)\|_2  + \|[D(0:t) \1]_i h(t)^\top \1 \|_2 }{n[D(0:t) \1]_i + n^2 \phi_i(t)}}\\
    & \;\;\;\; {\color{black} + \frac{n\sum_{l=1}^{t} \alpha(l) (\max_k [D(l:t)]_{ik}) \sum_{k=1}^n \| g_k(l)\|_2  }{n[D(0:t) \1]_i + n^2 \phi_i(t)}}\\
    & \le {\color{black}\frac{1}{n \eta} \Big[ n (\max_j | [D(0:t)]_{ij} |)  \sum_{k=1}^n \| h_k(0)\|_2} \\
    &\;\;\;\; {\color{black}+ n\sum_{l=1}^{t} \alpha(l) (\max_j | [D(l:t)]_{ij} |)   \sum_{k=1}^n \| g_k(l)\|_2 }\\
    &\;\;\;\; +n (\max_j | [D(0:t)]_{ij} |)\| h(t)^\top \1 \|_2\Big] \\
    & \le {\color{black}\frac{1}{ \eta} \Big[  4 \mu^{t} \sum_{k=1}^n  \|h_k(0)\|_2 + \sum_{l=1}^{t} \alpha(l) 4 \mu^{t-l} \sum_{k=1}^n  \|g_k(l)\|_2 }\\
    &\;\;\;\; + 4 \mu^{t} \| h(t)^\top \1 \|_2\Big]
\end{align*}
In addition, from \eqref{eq:h_update_1}, we have
$
    \| \1^\top h(t+1) \|_2
     \le \| \1^\top h(0)\|_2 + \| \sum_{l=1}^{t+1} \alpha(l) \1^\top g(l)\|_2 .
$
Then,
\begin{align*}
    &\;\;\;\;\; \| z_i(t+1) - \frac{ h(t)^\top \1}{n} \|_2 \\
    & \le 
    {\color{black}
    \frac{4}{ \eta} \Bigg[   \mu^{t}  \sum_{k=1}^n \| h_k(0)\|_2 + \sum_{l=1}^{t} \alpha(l)  \mu^{t-l} \sum_{k=1}^n  \|g_k(l)\|_2 }\\
    &\;\;\;\; +  \mu^{t} \| \1^\top h(0)\|_2 +  \mu^{t} \| \sum_{l=1}^{t} \alpha(l) \1^\top g(l)\|_2 \Bigg]\\
    & \le {\color{black}\frac{4}{ \eta} \left[   2 \mu^{t}  \sum_{k=1}^n \| h_k(0)\|_2 + 2 \sum_{l=0}^{t} \alpha(l)  \mu^{t-l} \sum_{k=1}^n \| g_k(l)\|_2  \right] }\\
    & = {\color{black}\frac{8}{ \eta} \Bigg[ \mu^{t} \sum_{k=1}^n \| x_k(0) - \alpha(0) g_k(0) \|_2 }\\
    &\;\;\;\; {\color{black} +  \sum_{l=0}^{t} \alpha(l)  \mu^{t-l} \sum_{k=1}^n \| g_k(l)\|_2  \Bigg] }
\end{align*}
Then, for all $t \ge 0$ and $i \in \mathcal{V}$,
\begin{align*}
    & \;\;\;\;\; \| z_i(t+1) - \frac{1}{n} \sum_{j=1}^n (x_j(t) -\alpha(t) g_j(t))\|_2 \\
    &\le
    {\color{black} \frac{8}{\eta} \Bigg[ \mu^t  \sum_{i=1}^n \| x_i(0) - \alpha(0) g_i(0) \|_2 }\\
    &\;\;\;\; 
    {\color{black}+ \sum_{s=0}^t  \mu^{t-s} \alpha(s) \sum_{i=1}^n\| g_i(s)\|_2 \Bigg] }\\
    &{\color{black} \le \frac{8}{\eta}  \mu^t  \sum_{i=1}^n \| x_i(0) - \alpha(0) g_i(0) \|_2 + \frac{8nG}{\eta} \sum_{s=0}^t  \mu^{t-s} \alpha(s)  }.
\end{align*}
If the stepsize sequence $\{ \alpha(t) \}$ satisfies Assumption~\ref{assum:step-size},
\begin{align*}
    & \;\;\;\;\; \| z_i(t+1) - \frac{1}{n} \sum_{j=1}^n (x_j(t) -\alpha(t) g_j(t))\|_2 \\
    &{\color{black} \le \frac{8}{\eta}  \mu^t  \sum_{i=1}^n \| x_i(0) - \alpha(0) g_i(0) \|_2 }\\
    &\;\;\;\;+  \frac{8nG}{\eta} \left(\sum_{s=0}^{\floor{\frac{t}{2}}}  \mu^{t-s} \alpha(s)
    + \sum_{s=\ceil{\frac{t}{2}}}^{t}  \mu^{t-s} \alpha(s) \right)   \\
    & {\color{black}\le \frac{8}{\eta}  \mu^t  \sum_{i=1}^n \| x_i(0) - \alpha(0) g_i(0) \|_2 }\\
    &\;\;\;\; +  \frac{8nG}{\eta(1-\mu)}   \left(  \alpha(0) \mu^{t/2}  +  \alpha(\ceil{\frac{t}{2}}) \right).
\end{align*}
This completes the proof.
\hfill $\qed$

\vspace{.05in}

We are now in a position to prove Theorem~\ref{thm:bound_everage_n_convex_bound}.

\vspace{.05in}

{\em Proof of Theorem~\ref{thm:bound_everage_n_convex_bound}:}
From Lemma~\ref{lemma:bound_consensus_push_SA}, for all $t \ge 0$ and $i \in \mathcal{V}$, 
\begin{align}
    &\;\;\;\;\; \|\langle z(t+1) \rangle - z_i(t+1) \|_2+\|\bar z(t+1) - z_i(t+1) \|_2 \nonumber\\ 
    & \le \|\langle z(t+1) \rangle - \frac{1}{n} \sum_{k=1}^n (x_k(t) -\alpha(t) g_k(t)) \|_2 \nonumber\\
    &\;\;\;\;+ \|\bar z(t+1) - \frac{1}{n} \sum_{k=1}^n (x_k(t) -\alpha(t) g_k(t)) \|_2 \nonumber\\
    &\;\;\; + 2 \| z_i(t+1) - \frac{1}{n} \sum_{k=1}^n (x_k(t) -\alpha(t) g_k(t)) \|_2, \nonumber
\end{align}
which implies that
\begin{align}
    &\;\;\;\;\; \|\langle z(t+1) \rangle - z_i(t+1) \|_2+\|\bar z(t+1) - z_i(t+1) \|_2 \nonumber\\ 
    & \le {\color{black} \sum_{j=1}^n (\pi_j(t) + \frac{1}{n} ) \|  z_j (t+1) - \frac{1}{n} \sum_{k=1}^n (x_k(t) -\alpha(t) g_k(t)) \|_2 }\nonumber\\
    &\;\;\; {\color{black} + 2 \| z_i(t+1) - \frac{1}{n} \sum_{k=1}^n (x_k(t) -\alpha(t) g_k(t)) \|_2,} \nonumber\\
    & {\color{black} \le \frac{32}{\eta}  \mu^t  \sum_{i=1}^n \| x_i(0) - \alpha(0) g_i(0) \|_2 +  \frac{32nG}{\eta}   \sum_{s=0}^t  \mu^{t-s} \alpha(s)}. \label{eq:bound_for_difference2_general}
\end{align}
In addition, when the stepsize sequence $\{ \alpha(t) \}$ satisfies Assumption~\ref{assum:step-size},
\begin{align}
    &\;\;\;\;\; \|\langle z(t+1) \rangle - z_i(t+1) \|_2+\|\bar z(t+1) - z_i(t+1) \|_2 \nonumber\\ 
    & {\color{black} \le \frac{32}{\eta}  \mu^t  \sum_{i=1}^n \| x_i(0) - \alpha(0) g_i(0) \|_2} \nonumber\\
    &\;\;\;\; +  \frac{32nG}{\eta(1-\mu)}   \left(  \alpha(0) \mu^{t/2}  -  \alpha(\ceil{\frac{t}{2}}) \right). \label{eq:bound_for_difference2}
\end{align}
From the update of $\langle z(t) \rangle$, we have
\begin{align}
    &\;\;\;\;\; \| \langle z(t+1) \rangle - z^* \|_2^2= \| \langle z(t) \rangle - z^* - \frac{\alpha(t)}{n}\sum_{i=1}^n g_i(t) \|_2^2 \nonumber\\
    & \le \| \langle z(t) \rangle - z^* \|_2^2 + \| \frac{\alpha(t)}{n}\sum_{i=1}^n g_i(t) \|_2^2 \nonumber\\
    &\;\;\;\; - 2 (\langle z(t) \rangle - z^*)^\top (\frac{\alpha(t)}{n}\sum_{i=1}^n g_i(t)) \nonumber \\
    & \le \| \langle z(t) \rangle - z^* \|_2^2 + \alpha^2(t) G^2 \nonumber\\
    &\;\;\;\; - 2 (\langle z(t) \rangle - z^*)^\top (\frac{\alpha(t)}{n}\sum_{i=1}^n g_i(t)). \label{eq:equation 1}
\end{align}
In addition,
\begin{align}
    &\;\;\;\;\; (\langle z(t) \rangle - z^*)^\top g_i(t)\nonumber\\
    & = (\langle z(t) \rangle - z_i(t) + z_i(t) - z^*)^\top g_i(t) \nonumber \\
    & = (\langle z(t) \rangle - z_i(t) )^\top g_i(t) + ( z_i(t) - z^*)^\top g_i(t) \nonumber \\
    & \ge f_i(z_i(t))  - f_i(z^*) - G \|\langle z(t) \rangle - z_i(t) \|_2  \label{eq:equation 4_1} \\ 
    & \ge f_i(\bar z(t))  - f_i(z^*) - G \|\langle z(t) \rangle - z_i(t) \|_2 \nonumber\\
    &\;\;\;\; - G \|\bar z(t) - z_i(t) \|_2, \label{eq:equation 4}
\end{align}
where we used \eqref{eq:subgradient} and \eqref{eq:G} in deriving \eqref{eq:equation 4_1}, and made use of \eqref{eq:G} to get \eqref{eq:equation 4}.

Combining \eqref{eq:equation 1} and \eqref{eq:equation 4}, we have
\begin{align*}
    &\;\;\;\;\; \| \langle z(t+1) \rangle - z^* \|_2^2  \\
    & \le \| \langle z(t) \rangle - z^* \|_2^2 + \alpha^2(t) G^2 - 2\alpha(t) ( f(\bar z(t) )  - f(z^*) )\\
    & \;\;\; + \frac{2G\alpha(t)}{n} \sum_{i=1}^n \left( \|\langle z(t) \rangle - z_i(t) \|_2 + \|\bar z(t) - z_i(t) \|_2 \right),
\end{align*}
which implies that
\begin{align*}
    &\;\;\;\;\;  2\alpha(t) ( f(\bar z(t) )  - f(z^*) )\\
    & \le \| \langle z(t) \rangle - z^* \|_2^2 + \alpha^2(t) G^2 - \| \langle z(t+1) \rangle - z^* \|_2^2 \\
    & \;\;\; + \frac{2G\alpha(t)}{n} \sum_{i=1}^n \left( \|\langle z(t) \rangle - z_i(t) \|_2 + \|\bar z(t) - z_i(t) \|_2 \right).
\end{align*}
Summing this up, we obtain
\begin{align*}
    &\;\;\;\;\; \sum_{\tau =0}^t 2\alpha(\tau) ( f(\bar z(\tau) )  - f(z^*) )\\
    & \le \| \langle z(0) \rangle - z^* \|_2^2  - \| \langle z(t+1) \rangle - z^* \|_2^2 + \sum_{\tau =0}^t \alpha^2(\tau) G^2\\
    & \;\;\; + \sum_{\tau =0}^t  \frac{2G\alpha(\tau)}{n} \sum_{i=1}^n \left( \|\langle z(\tau) \rangle - z_i(\tau) \|_2 + \|\bar z(\tau) - z_i(\tau) \|_2 \right).
\end{align*}
In addition, since
\begin{align*}
    &\;\;\;\; f\left(\frac{\sum_{\tau =0}^t \alpha(\tau) \bar z(\tau) }{\sum_{\tau =0}^t \alpha(\tau)}\right) - f(z^*)\nonumber\\
    & \le \frac{ \sum_{\tau =0}^t 2\alpha(\tau) ( f(\bar z(\tau)) - f(z^*) ) }{\sum_{\tau =0}^t 2\alpha(\tau)},
\end{align*}
then 
\begin{align}
    &\;\;\;\; f\left(\frac{\sum_{\tau =0}^t \alpha(\tau) \bar z(\tau) }{\sum_{\tau =0}^t \alpha(\tau)}\right) - f(z^*)\nonumber\\
    & \le \frac{ \| \langle z(0) \rangle - z^* \|_2^2  - \| \langle z(t+1) \rangle - z^* \|_2^2 + \sum_{\tau =0}^t \alpha^2(\tau) G^2}{\sum_{\tau =0}^t 2\alpha(\tau) }  \nonumber\\
    & + \frac{\sum_{\tau =0}^t \frac{2G\alpha(\tau)}{n} \sum_{i=1}^n (\|\langle z(\tau) \rangle - z_i(\tau) \|_2 + \|\bar z(\tau)  - z_i(\tau) \|_2)}{\sum_{\tau =0}^t 2\alpha(\tau) } \nonumber\\
    & \le  \frac{\sum_{\tau =0}^t {G\alpha(\tau)} \sum_{i=1}^n (\|\langle z(\tau) \rangle - z_i(\tau) \|_2 + \|\bar z(\tau)  - z_i(\tau) \|_2)}{{n}\sum_{\tau =0}^t \alpha(\tau) }  \nonumber\\
    & +  \frac{ \| \langle z(0) \rangle - z^* \|_2^2 + \sum_{\tau =0}^t \alpha^2(\tau) G^2}{\sum_{\tau =0}^t 2\alpha(\tau) } \label{eq:bound_mid}.
\end{align} 
We next consider the time-varying and fixed stepsizes separately. 

1) If the stepsize $\alpha(t)$ is time-varying and satisfies Assumption~\ref{assum:step-size}, then combining \eqref{eq:bound_for_difference2} and \eqref{eq:bound_mid}, we have
\begin{align*}
    &\;\;\;\;\; f\left(\frac{\sum_{\tau =0}^t \alpha(\tau) \bar z(\tau) }{\sum_{\tau =0}^t \alpha(\tau)}\right) - f(z^*) \\
    & \le \frac{ \| \langle z(0) \rangle - z^* \|_2^2 + G^2\sum_{\tau =0}^t \alpha^2(\tau) }{\sum_{\tau =0}^t 2\alpha(\tau) }  \\
    &\;\;\;\;  + \frac{ G\alpha(0) \sum_{i=1}^n (\|\langle z(0) \rangle - z_i(0) \|_2 + \|\bar z(0)  - z_i(0) \|_2)}{n \sum_{\tau =0}^t \alpha(\tau) } \\
     &\;\;\;\; {\color{black}+ \frac{32G}{\eta}  \sum_{i=1}^n \| x_i(0) - \alpha(0) g_i(0) \|_2 \frac{\sum_{\tau =0}^{t-1} \alpha(\tau) \mu^\tau  }{\sum_{\tau =0}^t \alpha(\tau) } }\\
     &\;\;\;\; + \frac{32nG^2}{\eta(1-\mu)} \frac{\sum_{\tau =0}^{t-1} \alpha(\tau)   \left(  \alpha(0) \mu^{\tau/2}  +  \alpha(\ceil{\frac{\tau}{2}}) \right)}{\sum_{\tau =0}^t \alpha(\tau) }.
\end{align*} 
For all $i\in\mathcal{V}$, from Proposition~\ref{lemma:push-sum_pi_intfty} and $y_i(0) = 1$, we have $\pi_i(0) = \frac{1}{n}$, which implies that $ \langle z(0) \rangle = \frac{1}{n}\sum_i^n z_i(0) = \bar z(0)$. We thus have derived \eqref{eq:bound_timevarying}.

2) If the stepsize is fixed and  $\alpha(t) = 1/\sqrt{T}$ for all $t\ge0$, then from \eqref{eq:bound_mid}, we have
\begin{align*}
    &\;\;\;\;\; f\left(\frac{\sum_{\tau =0}^{T-1} \bar z(\tau)  }{ T } \right) - f(z^*)\\
    & \le \frac{{G} \sum_{\tau =0}^{T-1} \sum_{i=1}^n \|\langle z(\tau) \rangle - z_i(\tau) \|_2+\|\bar z(\tau) - z_i(\tau) \|_2}{{n}T }\\
    &\;\;\;\; + \frac{ \| \langle z(0) \rangle - z^* \|_2^2 + G^2 }{ 2\sqrt{T} }.
\end{align*}
Using \eqref{eq:bound_for_difference2_general}, we have
\begin{align*}
    &\;\;\;\; f\left(\frac{\sum_{\tau =0}^{T-1} \bar z(\tau)  }{ T } \right) - f(z^*) \\
    & \le \frac{G \sum_{i=1}^n \|\langle z(0) \rangle - z_i(0) \|_2+\|\bar z(0) - z_i(0) \|_2}{nT }\\ 
    &\;\;\;\; + \frac{ \| \langle z(0) \rangle - z^* \|_2^2 + G^2 }{ 2\sqrt{T} }+ \frac{32n G^2 }{T \eta} \sum_{\tau =0}^{T-2} \sum_{s=0}^\tau  \mu^{\tau-s} \frac{1}{\sqrt{T}}\\
    &\;\;\;\; {\color{black}+ \frac{32 G }{T \eta}  \sum_{i=1}^n \| x_i(0) - \frac{1}{\sqrt{T}} g_i(0) \|_2 \sum_{\tau =0}^{T-2} \mu^\tau ,}
\end{align*}
which implies that
\begin{align*}
    &\;\;\;\; f\left(\frac{\sum_{\tau =0}^{T-1} \bar z(\tau)  }{ T } \right) - f(z^*) \\
    & \le  \frac{G \sum_{i=1}^n \|\langle z(0) \rangle - z_i(0) \|_2+\|\bar z(0) - z_i(0) \|_2}{{n}T }\\
    &\;\;\;\; +\frac{ \| \langle z(0) \rangle - z^* \|_2^2 + G^2 }{ 2\sqrt{T} }+ \frac{32n G^2 }{ \sqrt{T}\eta(1- \mu)} \\
    &\;\;\; {\color{black}+ \frac{32 G }{T \eta(1 - \mu)} \sum_{i=1}^n \|  x_i(0) - \frac{1}{\sqrt{T}} g_i(0) \|_2 .}
\end{align*}
Since $ \langle z(0) \rangle = \frac{1}{n}\sum_i^n z_i(0) = \bar z(0)$, we have derived \eqref{eq:bound_fixed}.
This completes the proof.
\hfill $\qed$

\vspace{.05in}

We next prove Theorem~\ref{thm:bound_everage_n_convex}.

\vspace{.05in}

{\em Proof of Theorem~\ref{thm:bound_everage_n_convex}:}
1) If the stepsize $\alpha(t)$ is time-varying and satisfies Assumption~\ref{assum:step-size}, then 
\begin{align*}
    & \lim_{t\to\infty}\frac{ \| \langle z(0) \rangle - z^* \|_2^2 + \sum_{\tau =0}^t \alpha^2(\tau) G^2}{\sum_{\tau =0}^t 2\alpha(\tau) } = 0\\
    &\lim_{t\to\infty} \frac{ \sum_{i=1}^n (\|\langle z(0) \rangle - z_i(0) \|_2 + \|\bar z(0)  - z_i(0) \|_2)}{\sum_{\tau =0}^t \alpha(\tau) } = 0.
\end{align*}
In addition, since $\sum_{\tau =0}^{t-1} \alpha(\tau) \mu^\tau \le \frac{\alpha(0)   }{1-\mu}$ and 
\begin{align*}
    &\;\;\;\;\; \sum_{\tau =0}^{t-1} \alpha(\tau)   \left(  \alpha(0) \mu^{\tau/2}  +  \alpha(\ceil{\frac{\tau}{2}}) \right) \\
    &\le \alpha(0)^2 \sum_{\tau =0}^{t-1} \mu^{\tau/2}  +  \sum_{\tau =0}^{t-1} \alpha(\ceil{\frac{\tau}{2}})^2 \\
    &\le \frac{\alpha(0)^2 }{1- \mu^{1/2}}  +  \sum_{\tau =0}^{t-1} \alpha(\ceil{\frac{\tau}{2}})^2,
\end{align*}
it follows that
\begin{align*}
    & \lim_{t\to\infty} \frac{\sum_{\tau =0}^{t-1} \alpha(\tau) \mu^\tau  }{\sum_{\tau =0}^t \alpha(\tau) } = 0\\
    &\lim_{t\to\infty} \frac{\sum_{\tau =0}^{t-1} \alpha(\tau)   \left(  \alpha(0) \mu^{\tau/2}  +  \alpha(\ceil{\frac{\tau}{2}}) \right)}{\sum_{\tau =0}^t \alpha(\tau) } = 0.
\end{align*}
From \eqref{eq:bound_timevarying}, we have 
\begin{align*}
    \lim_{t\to\infty} f\left(\frac{\sum_{\tau =0}^t \alpha(\tau) \bar z(\tau) }{\sum_{\tau =0}^t \alpha(\tau)}\right) - f(z^*) = 0.
\end{align*}
2) If the stepsize is fixed and  $\alpha(t) = 1/\sqrt{T}$ for all $t\ge0$, then from \eqref{eq:bound_fixed}, we have
\begin{align*}
    f\bigg(\frac{\sum_{\tau =0}^{T-1} \bar z(\tau)  }{ T } \bigg) - f(z^*) 
    & \le O\Big(\frac{1}{ \sqrt{T}}\Big).
\end{align*}
This completes the proof.
\hfill $\qed$

\vspace{.05in}

We finally prove Theorem~\ref{thm:bound_everage_zi}.

\vspace{.05in}

{\em Proof of Theorem~\ref{thm:bound_everage_zi}:}
From Lemma~\ref{lemma:bound_consensus_push_SA}, for all $t \ge 0$ and $i,j \in \mathcal{V}$,
\begin{align}
    &\;\;\;\;\; \|\langle z(t+1) \rangle - z_i(t+1) \|_2+\| z_j(t+1) - z_i(t+1) \|_2 \nonumber\\ 
    & \le \|\langle z(t+1) \rangle - \frac{1}{n} \sum_{k=1}^n (x_k(t) -\alpha(t) g_k(t)) \|_2 \nonumber\\
    &\;\;\;\;+ \|z_j(t+1) - \frac{1}{n} \sum_{k=1}^n (x_k(t) -\alpha(t) g_k(t)) \|_2 \nonumber\\
    &\;\;\;\; + 2 \| z_i(t+1) - \frac{1}{n} \sum_{k=1}^n (x_k(t) -\alpha(t) g_k(t)) \|_2, \nonumber\\
    &{\color{black} \le \frac{32}{\eta}  \mu^t  \sum_{i=1}^n \| x_i(0) - \alpha(0) g_i(0) \|_2 +  \frac{32nG}{\eta}   \sum_{s=0}^t  \mu^{t-s} \alpha(s).} \label{eq:bound_for_difference2_general_zi}
\end{align}
When the stepsize sequence $\{ \alpha(t) \}$ satisfies Assumption~\ref{assum:step-size},
\begin{align}
    &\;\;\;\;\; \|\langle z(t+1) \rangle - z_i(t+1) \|_2+\|z_j(t+1) - z_i(t+1) \|_2 \nonumber\\ 
    &{\color{black}\le \frac{32}{\eta}  \mu^t \sum_{i=1}^n \|  x_i(0) - \alpha(0) g_i(0) \|_2} \nonumber\\
    &\;\;\;\; +  \frac{32nG}{\eta(1-\mu)}   \left(  \alpha(0) \mu^{t/2}  +  \alpha(\ceil{\frac{t}{2}}) \right). \label{eq:bound_for_difference2_zi}
\end{align}
From \eqref{eq:subgradient}, \eqref{eq:G} and \eqref{eq:equation 4_1},
we have for any $k\in \scr V$
\begin{align}
    &\;\;\;\;\; (\langle z(t) \rangle - z^*)^\top g_i(t)\nonumber\\
    & \ge f_i(z_i(t))  - f_i(z^*) - G \|\langle z(t) \rangle - z_i(t) \|_2  \nonumber \\ 
    & \ge f_i( z_k(t))  - f_i(z^*) - G \|\langle z(t) \rangle - z_i(t) \|_2 \nonumber\\
    &\;\;\;\; - G \| z_k(t) - z_i(t) \|_2.
    \label{eq:equation 4_zi}
\end{align}
Then, combining \eqref{eq:equation 1} and \eqref{eq:equation 4_zi}, we have
\begin{align*}
    &\;\;\;\;\; \| \langle z(t+1) \rangle - z^* \|_2^2  \\
    & \le \| \langle z(t) \rangle - z^* \|_2^2 + \alpha^2(t) G^2 - 2\alpha(t) ( f(z_k(t) )  - f(z^*) )\\
    & \;\;\; + \frac{2G\alpha(t)}{n} \sum_{i=1}^n \left( \|\langle z(t) \rangle - z_i(t) \|_2 + \| z_k(t) - z_i(t) \|_2 \right),
\end{align*}
which implies that
\begin{align*}
    &\;\;\;\;\;  2\alpha(t) ( f(z_k(t) )  - f(z^*) )\\
    & \le \| \langle z(t) \rangle - z^* \|_2^2 + \alpha^2(t) G^2 - \| \langle z(t+1) \rangle - z^* \|_2^2 \\
    & \;\;\; + \frac{2G\alpha(t)}{n} \sum_{i=1}^n \left( \|\langle z(t) \rangle - z_i(t) \|_2 + \|z_k(t) - z_i(t) \|_2 \right). 
\end{align*}
Summing this up, we obtain
\begin{align*}
    &\;\;\;\;\; \sum_{\tau =0}^t 2\alpha(\tau) ( f(z_k(\tau) )  - f(z^*) )\\
    & \le \| \langle z(0) \rangle - z^* \|_2^2  - \| \langle z(t+1) \rangle - z^* \|_2^2 + \sum_{\tau =0}^t \alpha^2(\tau) G^2\\
    & \;\;\; + \sum_{\tau =0}^t  \frac{2G\alpha(\tau)}{n} \sum_{i=1}^n \left( \|\langle z(\tau) \rangle - z_i(\tau) \|_2 + \| z_k(\tau) - z_i(\tau) \|_2 \right).
\end{align*}
In addition, 
\begin{align}
    &\;\;\;\; f\left(\frac{\sum_{\tau =0}^t \alpha(\tau) z_k(\tau) }{\sum_{\tau =0}^t \alpha(\tau)}\right) - f(z^*)\nonumber\\
    & \le \frac{ \sum_{\tau =0}^t 2\alpha(\tau) ( f(z_k(\tau)) - f(z^*) ) }{\sum_{\tau =0}^t 2\alpha(\tau)}  \nonumber\\
    & \le \frac{ \| \langle z(0) \rangle - z^* \|_2^2  - \| \langle z(t+1) \rangle - z^* \|_2^2 + \sum_{\tau =0}^t \alpha^2(\tau) G^2}{\sum_{\tau =0}^t 2\alpha(\tau) }  \nonumber\\
    & + \frac{\sum_{\tau =0}^t \frac{2G\alpha(\tau)}{n} \sum_{i=1}^n (\|\langle z(\tau) \rangle - z_i(\tau) \|_2 + \|z_k(\tau)  - z_i(\tau) \|_2)}{\sum_{\tau =0}^t 2\alpha(\tau) } \nonumber\\
    & \le  \frac{\sum_{\tau =0}^t G\alpha(\tau) \sum_{i=1}^n (\|\langle z(\tau) \rangle - z_i(\tau) \|_2 + \|z_k(\tau)  - z_i(\tau) \|_2)}{{n}\sum_{\tau =0}^t \alpha(\tau) }  \nonumber\\
    & +  \frac{ \| \langle z(0) \rangle - z^* \|_2^2 + \sum_{\tau =0}^t \alpha^2(\tau) G^2}{\sum_{\tau =0}^t 2\alpha(\tau) } \label{eq:bound_mid_zi}.
\end{align} 

1) If the stepsize $\alpha(t)$ is time-varying and satisfies Assumption~\ref{assum:step-size}, then combining \eqref{eq:bound_for_difference2_zi} and \eqref{eq:bound_mid_zi}, we have
\begin{align*}
    &\;\;\;\;\; f\left(\frac{\sum_{\tau =0}^t \alpha(\tau) z_k(\tau) }{\sum_{\tau =0}^t \alpha(\tau)}\right) - f(z^*) \\
    & \le \frac{ \| \langle z(0) \rangle - z^* \|_2^2 + G^2\sum_{\tau =0}^t \alpha^2(\tau) }{\sum_{\tau =0}^t 2\alpha(\tau) }  \\
    &\;\;\;\;  + \frac{ G\alpha(0) \sum_{i=1}^n (\|\langle z(0) \rangle - z_i(0) \|_2 + \|z_k(0)  - z_i(0) \|_2)}{n\sum_{\tau =0}^t \alpha(\tau) } \\
     &\;\;\;\; {\color{black}+ \frac{32G}{\eta}  \sum_{i=1}^n \| x_i(0) - \alpha(0) g_i(0) \|_2 \frac{\sum_{\tau =0}^{t-1} \alpha(\tau) \mu^\tau  }{\sum_{\tau =0}^t \alpha(\tau) } }\\
     &\;\;\;\; + \frac{32nG^2}{\eta(1-\mu)} \frac{\sum_{\tau =0}^{t-1} \alpha(\tau)   \left(  \alpha(0) \mu^{\tau/2}  +  \alpha(\ceil{\frac{\tau}{2}}) \right)}{\sum_{\tau =0}^t \alpha(\tau) }.
\end{align*} 
Since$ \langle z(0) \rangle = \bar z(0)$, we have derived \eqref{eq:bound_timevarying_zi}.

2) If the stepsize is fixed and  $\alpha(t) = 1/\sqrt{T}$ for all $t\ge0$, then from \eqref{eq:bound_mid_zi}, we have
\begin{align*}
    &\;\;\;\;\; f\left(\frac{\sum_{\tau =0}^{T-1} z_k(\tau)  }{ T } \right) - f(z^*)\\
    & \le \frac{ G \sum_{\tau =0}^{T-1} \sum_{i=1}^n \|\langle z(\tau) \rangle - z_i(\tau) \|_2+\|z_k(\tau) - z_i(\tau) \|_2}{nT }\\
    &\;\;\;\; + \frac{ \| \langle z(0) \rangle - z^* \|_2^2 + G^2 }{ 2\sqrt{T} }.
\end{align*}
Using \eqref{eq:bound_for_difference2_general_zi}, we have
\begin{align*}
    &\;\;\;\; f\left(\frac{\sum_{\tau =0}^{T-1} z_k(\tau)  }{ T } \right) - f(z^*) \\
    & \le \frac{G \sum_{i=1}^n \|\langle z(0) \rangle - z_i(0) \|_2+\|z_k(0) - z_i(0) \|_2}{nT }\\ 
    &\;\;\;\; + \frac{ \| \langle z(0) \rangle - z^* \|_2^2 + G^2 }{ 2\sqrt{T} }+ \frac{32n G^2 }{T \eta} \sum_{\tau =0}^{T-2} \sum_{s=0}^\tau  \mu^{\tau-s} \frac{1}{\sqrt{T}}\\
    &\;\;\;\; {\color{black}+ \frac{32 G }{T \eta}  \sum_{i=1}^n \| x_i(0) - \frac{1}{\sqrt{T}} g_i(0) \|_2 \sum_{\tau =0}^{T-2} \mu^\tau,}
\end{align*}
which implies that
\begin{align*}
    &\;\;\;\; f\left(\frac{\sum_{\tau =0}^{T-1} z_k(\tau)  }{ T } \right) - f(z^*) \\
    & \le  \frac{G \sum_{i=1}^n \|\langle z(0) \rangle - z_i(0) \|_2+\|z_k(0) - z_i(0) \|_2}{nT }\\
    &\;\;\;\; +\frac{ \| \langle z(0) \rangle - z^* \|_2^2 + G^2 }{ 2\sqrt{T} }+ \frac{32n G^2 }{ \sqrt{T}\eta(1- \mu)} \\
    &\;\;\; {\color{black}+ \frac{32 G }{T \eta(1 - \mu)}  \sum_{i=1}^n \| x_i(0) - \frac{1}{\sqrt{T}} g_i(0) \|_2 .}
\end{align*}
Since $ \langle z(0) \rangle = \bar z(0)$, we have derived  \eqref{eq:bound_fixed_zi}.
This completes the proof.
\hfill $\qed$
\vspace{.1in}

\section{Conclusion}

The well-know push-sum based subgradient algorithm for distributed convex optimization over unbalanced directed graphs has been revisited. A novel analysis tool has been proposed, which improves the convergence rate of the {\color{black} subgradient-push} algorithm from $O(\ln t/\sqrt{t})$ to   $O(1/\sqrt{t})$, which is the same as that of the single-agent subgradient method and thus optimal.
As a future work, the proposed tool is expected to be applicable to analyze other push-sum based algorithms and improve/simplify their convergence analyses, for example, DEXTRA \cite{xi2017dextra} and Push-DIGing \cite{nedic2017achieving}. 
Another future direction is to extend the proposal tool to push-sum based distributed algorithms with communication delays and asynchronous updating.

\vspace{.1in}

\bibliographystyle{unsrt}
\bibliography{push}

\begin{thebibliography}{10}

\bibitem{reachingp1}
M.~Cao, A.S. Morse, and B.D.O. Anderson.
\newblock Reaching a consensus in a dynamically changing environment: A
  graphical approach.
\newblock {\em SIAM Journal on Control and Optimization}, 47(2):575--600, 2008.

\bibitem{fast}
L.~Xiao and S.~Boyd.
\newblock Fast linear iterations for distributed averaging.
\newblock {\em Systems \& Control Letters}, 53(1):65--78, 2004.

\bibitem{push}
D.~Kempe, A.~Dobra, and J.~Gehrke.
\newblock Gossip-based computation of aggregate information.
\newblock In {\em Proceedings of the 44th IEEE Symposium on Foundations of
  Computer Science}, pages 482--491, 2003.

\bibitem{nedic}
A.~Nedi\'{c} and A.~Olshevsky.
\newblock Distributed optimization over time-varying directed graphs.
\newblock {\em IEEE Transactions on Automatic Control}, 60(3):601--615, 2015.

\bibitem{yixuan}
Y.~Lin, K.~Zhang, Z.~Yang, Z.~Wang, T.~Ba\c{s}ar, R.~Sandhu, and J.~Liu.
\newblock A communication-efficient multi-agent actor-critic algorithm for
  distributed reinforcement learning.
\newblock In {\em Proceedings of the 58th IEEE Conference on Decision and
  Control}, pages 5562--5567, 2019.

\bibitem{weighted}
F.~B\'{e}n\'{e}zit, V.~Blondel, P.~Thiran, J.~N. Tsitsiklis, and M.~Vetterli.
\newblock Weighted gossip: distributed averaging using non-doubly stochastic
  matrices.
\newblock In {\em Proceedings of the 2010 IEEE International Symposium on
  Information Theory}, pages 1753--1757, 2010.

\bibitem{hadjicostis2013average}
C.N. Hadjicostis and T.~Charalambous.
\newblock Average consensus in the presence of delays in directed graph
  topologies.
\newblock {\em IEEE Transactions on Automatic Control}, 59(3):763--768, 2013.

\bibitem{acc12}
J.~Liu and A.S. Morse.
\newblock Asynchronous distributed averaging using double linear iterations.
\newblock In {\em Proceedings of the 2012 American Control Conference}, pages
  6620--6625, 2012.

\bibitem{xi2017dextra}
C.~Xi and U.A. Khan.
\newblock {DEXTRA}: A fast algorithm for optimization over directed graphs.
\newblock {\em IEEE Transactions on Automatic Control}, 62(10):4980--4993,
  2017.

\bibitem{shi2015extra}
W.~Shi, Q.~Ling, G.~Wu, and W.~Yin.
\newblock {EXTRA}: An exact first-order algorithm for decentralized consensus
  optimization.
\newblock {\em SIAM Journal on Optimization}, 25(2):944--966, 2015.

\bibitem{nedic2017achieving}
A.~Nedi\'c, A.~Olshevsky, and W.~Shi.
\newblock Achieving geometric convergence for distributed optimization over
  time-varying graphs.
\newblock {\em SIAM Journal on Optimization}, 27(4):2597--2633, 2017.

\bibitem{nedic2009distributed}
A.~Nedi\'c and A.~Ozdaglar.
\newblock Distributed subgradient methods for multi-agent optimization.
\newblock {\em IEEE Transactions on Automatic Control}, 54(1):48--61, 2009.

\bibitem{yang2019survey}
T.~Yang, X.~Yi, J.~Wu, Y.~Yuan, D.~Wu, Z.~Meng, Y.~Hong, H.~Wang, Z.~Lin, and
  K.H Johansson.
\newblock A survey of distributed optimization.
\newblock {\em Annual Reviews in Control}, 47:278--305, 2019.

\bibitem{nedic2018distributed}
A.~Nedi{\'c} and J.~Liu.
\newblock Distributed optimization for control.
\newblock {\em Annual Review of Control, Robotics, and Autonomous Systems},
  1:77--103, 2018.

\bibitem{molzahn2017survey}
D.K. Molzahn, F.~D{\"o}rfler, H.~Sandberg, S.H. Low, S.~Chakrabarti,
  R.~Baldick, and J.~Lavaei.
\newblock A survey of distributed optimization and control algorithms for
  electric power systems.
\newblock {\em IEEE Transactions on Smart Grid}, 8(6):2941--2962, 2017.

\bibitem{gharesifard2013distributed}
B.~Gharesifard and J.~Cort{\'e}s.
\newblock Distributed continuous-time convex optimization on weight-balanced
  digraphs.
\newblock {\em IEEE Transactions on Automatic Control}, 59(3):781--786, 2013.

\bibitem{nedic2018network}
A.~Nedi{\'c}, A.~Olshevsky, and M.~G Rabbat.
\newblock Network topology and communication-computation tradeoffs in
  decentralized optimization.
\newblock {\em Proceedings of the IEEE}, 106(5):953--976, 2018.

\bibitem{subgradient}
B.~Polyak.
\newblock A general method for solving extremum problems.
\newblock {\em Doklady Akademii Nauk}, 8(3):593--597, 1967.

\bibitem{nedic2009distributed_quan}
A.~Nedi\'c, A.~Olshevsky, A.~Ozdaglar, and J.N. Tsitsiklis.
\newblock On distributed averaging algorithms and quantization effects.
\newblock {\em IEEE Transactions on Automatic Control}, 54(11):2506--2517,
  2009.

\bibitem{iutzeler2013analysis}
F.~Iutzeler, P.~Ciblat, and W.~Hachem.
\newblock Analysis of sum-weight-like algorithms for averaging in wireless
  sensor networks.
\newblock {\em IEEE Transactions on Signal Processing}, 61(11):2802--2814,
  2013.

\bibitem{kolmogorov}
A.~Kolmogoroff.
\newblock Zur theorie der markoffschen ketten.
\newblock {\em Mathematische Annalen}, 112(1):155--160, 1936.

\bibitem{blackwell}
D.~Blackwell.
\newblock Finite non-homogeneous chains.
\newblock {\em Annals of Mathematics}, 46(4):594--599, 1945.

\bibitem{tacrate}
A.~Nedi\'c and J.~Liu.
\newblock On convergence rate of weighted-averaging dynamics for consensus
  problems.
\newblock {\em IEEE Transactions on Automatic Control}, 62(2):766--781, 2017.

\bibitem{touri2012product}
B.~Touri.
\newblock {\em Product of Random Stochastic Matrices and Distributed
  Averaging}.
\newblock Springer Science \& Business Media, 2012.

\end{thebibliography}

\end{document}